\documentclass[11pt]{article} 
\usepackage[latin1]{inputenc}
\usepackage[T1]{fontenc}               
\usepackage[francais]{babel}

\usepackage{amssymb}
\usepackage{amsmath}
\usepackage{mathrsfs}
\newtheorem{theo}{Th\'eor\`eme}[section]
\newtheorem{lemme}[theo]{Lemme}
\newtheorem{prop}[theo]{Proposition}
\newtheorem{cor}[theo]{Corollaire}
\newtheorem{rem}[theo]{Remarque }
\newtheorem{exe}[theo]{Exemple }
\newenvironment{remarque}{\begin{rem}\em}{\end{rem}}
\newenvironment{exemple}{\begin{exe}\em}{\end{exe}}
\newtheorem{rems}[theo]{Remarques }

\newtheorem{defi}[theo]{D\'efinition } 

\newenvironment{dem}{\noindent{\it D\'emonstration}.}
{\unskip\hfill\null\nobreak\hfill\carre\vskip1em\par}
\newcommand{\carre}{\rule{1ex}{1ex}}
\newcommand{\rond}{\raisebox{.3mm}{\scriptsize$\circ$}}

\def\g{\mathfrak}
\def\R{\mathbb R}
\def\C{\mathbb C}

\def\N{\mathbb N}

\def\P{\mathcal P}
\def\SL{\mathop{\mbox{\rm SL}}\nolimits}

\def\ad{\mathop{\mbox{\rm ad}}\nolimits}

\def\id{\mathop{\mbox{\rm id}}\nolimits}
\def\Car{\mathop{\mbox{\rm Car}}\nolimits}
\def\diag{\mathop{\mbox{\rm diag}}\nolimits}
\def\Der{\mathop{\mbox{\rm Der}}\nolimits}
\def\inv{\mathop{\scriptsize{\mbox{\rm inv}}}\nolimits}
\def\reg{\mathop{\scriptsize{\mbox{\rm reg}}}\nolimits}
\title{Sur les champs de vecteurs invariants  sur l'espace
tangent d'un espace sym\'etrique r\'eductif}
\author{Abderrazak Bouaziz\footnote{Universit\'e de Poitiers-CNRS, Laboratoire de Math\'ematiques et
Applications, BP 30179, 86962 Futuroscope-Chasseneuil, France, 
bouaziz@math.univ-poitiers.fr} \quad  et \quad Nouri
  Kamoun\footnote{Facult\'e des Sciences de Monastir, 5019 Monastir,
    Tunisie, nouri.kamoun@fsm.rnu.tn}}
\date{}

\begin{document}

\maketitle

\section{Introduction}
Soit $K$ un groupe de Lie op\'erant
diff\'erentiablement dans une vari\'et\'e diff\'erentiable $X$. Tout champ de vecteurs lisse et $K$-invariant
sur $X$ d\'efinit naturellement une d\'erivation de l'alg\`ebre
$C^{\infty}(X)^K$ des fonctions
lisses et $K$-invariantes sur $X$, mais on n'obtient pas ainsi toutes
les d\'erivations de cette alg\`ebre. Par exemple (\cite{Sch},
exemple 2.6) pour $X=\R$ et
$K=\{\pm\id_{\R}\}$, les fonctions invariantes par $K$ sont les
fonctions paires, et si $f\in C^{\infty}(\R)^{K}$, il existe
une fonction lisse $g$ sur $\R$ telle que $f(x)=g(x^2)$ pour tout
$x\in\R$ ; la d\'erivation $D$ de $C^{\infty}(\R)^{K}$
d\'efinie par $Df(x)=\frac{d}{dx}g (x^2)$ ne provient pas d'un champ de
vecteurs lisse $K$-invariant sur $\R$.

Lorsque $K$ est compact, G. Schwarz a donn\'e  une
caract\'erisation des d\'erivations de $C^{\infty}(X)^K$ qui proviennent de
champs de vecteurs $K$-invariants sur $X$ (\cite{Sch}, th\'eor\`eme 0.2) : ce sont les d\'erivations
qui pr\'eservent les id\'eaux de fonctions qui s'annulent sur les startes
de $X$ (si $L$ est un sous-groupe de $K$, l'ensemble des \'el\'ements de
$X$ dont le fixateur dans $K$ est conjugu\'e \`a $L$ est appel\'e une strate de
$X$). Dans (\cite{Sch},
remarque 6.16) il est dit qu'on pourrait  formuler et essayer de
d\'emontrer des r\'esultats analogues pour l'action de groupes r\'eductifs
r\'eels. Nous nous proposons d'\'etudier le cas des groupes d'isotropie des
espaces sym\'etriques r\'eductifs agissant sur les espaces tangents. Pour
cela, nous nous inspirons de la formulation de M. Ra\"is pour les
espaces sym\'etriques riemanniens.

Supposons que  $K$ (toujours compact) soit  le groupe d'isotropie d'un
espace sym\'etrique riemannien op\'erant (lin\'eairement) sur l'espace
tangent $\g p$ de cet espace sym\'etrique, et soient
$(P_1,\ldots,P_{\ell})$ un syst\`eme minimal de g\'en\'erateurs homog\`enes de
l'alg\`ebre des polyn\^omes $K$-invariants sur $\g p$. On note 
$\nabla P_i$ le gradient de $P_i$ relativement \`a la structure riemannienne $B$
de $\g p$, et on pose $\it{\Phi}=\det(B(\nabla p_i,\nabla p_j))$. 
Dans cette situation,  M. Ra\"is  
a montr\'e (\cite{Ra}, lemme 3.4) que les d\'erivations de $C^{\infty}(\g p)^K$ qui proviennent de
champs de  vecteurs invariants sont celles qui pr\'eservent
l'id\'eal engendr\'e par $\it{\Phi}$, et que toute  d\'erivation de ce type
est induite
par un champ de vecteurs de la forme
$\sum_{i=1}^{\ell}\varphi_i\nabla P_i$, avec $\varphi_i\in
C^{\infty}(\g p)^K$.

L'objectif de cet article est de montrer que la formulation   de
Ra\"is se g\'en\'eralise aux espaces sym\'etriques r\'eductifs
quelconques.  L'id\'ee de la preuve est la m\^eme que celle de Ra\"is,
cependant notre situation pr\'esente deux difficult\'es suppl\'ementaires :
l'existence de plusieurs classes de conjugaison de sous-espaces de
Cartan, et, surtout, l'ensemble des restrictions \`a un sous-espace de Cartan des
fonctions invariantes n'est pas l'ensemble des fonctions invariantes
par un groupe. Cela nous oblige \`a travailler avec les d\'eveloppements
de Taylor des fonctions invariantes qui eux pr\'esentent de bonnes
propri\'et\'es d'invariance  (\cite{Ka} et
\cite{Ka2}).

Dans le cas  sym\'etrique riemannien, L. Sliman-Bouattour a donn\'e  dans \cite{SB} une description de
l'ensemble $\mathscr Z$ des champs de vecteurs lisses et
$K$-invariants qui annulent tous les \'el\'ements de $C^{\infty}(\g p)^K$ : si l'espace sym\'etrique est irr\'eductible, il y a deux cas
qui se pr\'esentent, ou bien la repr\'esentation de $K$ dans le
complexifi\'e $\g p_{\C}$ de $\g p$ est irr\'eductible, dans ce cas
$\mathscr Z$ est r\'eduit \`a $\{0\}$, ou bien il existe une structure
complexe $K$-invariante $J:\g p\rightarrow\g p$, dans ce cas $\mathscr
Z$ est un $C^{\infty}(\g p)^K$-module libre dont $J\nabla
P_1,\ldots,J\nabla P_{\ell}$ est une base. En particulier $\mathscr Z$
est engendr\'e par des champs de vecteurs invariants \`a coefficients polynomiaux.
Nous verrons \`a la fin de cet article que pour la paire
$(\g{sl}(3,\R),\g{so}(2,1))$ l'ensemble des  \'el\'ements de  $\mathscr Z$ \`a coefficients polynomiaux n'engendre
pas  $\mathscr Z$ comme module sur les fonctions lisses invariantes.

Signalons enfin qu'on trouve dans \cite{A-G} une description des
champs de vecteurs r\'eguliers invariants sur les espaces sym\'etriques
r\'eductifs complexes, analogue \`a celle de Ra\"is et Sliman-Bouattour.

Nous remercions Mustapha Ra\"is pour ses commentaires sur une premi\`ere
version de cet article.

\section{Enonc\'e du  r\'esultat principal }\label{resultat}

Si $V$ est un espace vectoriel sur $\R$ ou $\C$, on note $\C[V]$
l'ensemble des  polyn\^omes
complexes sur $V$. Si $K$ est un sous-groupe de $GL(V)$, on note $V^K$
l'ensemble des \'el\'ements de $V$ fix\'es par tous les \'el\'ements de $K$. De
m\^eme si $\mathcal F$ est un espace de fonctions sur $V$, on note
$\mathcal F^K$ l'ensemble des \'el\'ements de $\mathcal F$ fix\'es par
l'action naturelle de $K$ dans $\mathcal F$.

Si $V$ est un espace vectoriel sur $\R$, on note $V_{\C}$ son
complexifi\'e. Les alg\`ebres $\C[V]$ et  $\C[V_{\C}]$ sont  
naturellement isomorphes, et seront identifi\'ees dans la suite. Ainsi
les notations $\C[V]$ et $\C[V_{\C}]$ d\'esigneront  la m\^eme chose ;
on utilisera la seconde pour souligner  l'action d'un
sous-groupe de $GL(V_{\C})$ dans les ployn\^omes. Aussi  
la m\^eme lettre d\'esignera \`a la fois un polyn\^ome sur
$V$ et son prolongement \`a $V_{\C}$.

Si $U$ est un ouvert de $V$, on note $C^{\infty}(U)$ l'ensembles des
fonctions complexes  de classe $C^{\infty}$ sur $U$. On note $\g
X(U)=C^{\infty}(U,V_{\C})$ l'ensemble des champs de vecteurs sur $U$ \`a
valeurs dans $V_{\C}$ ; ces champs de vecteurs op\`erent naturellement
dans $C^{\infty}(U)$ : si $v\in V_{\C}$, vu comme champ de vecteurs
constant, et si on \'ecrit $v=v_1+iv_2, v_i\in V$, alors, pour tout $f\in C^\infty(U)$,
$$v\cdot f(x)=\partial(v_1)f(x)+i\partial(v_2)f(x)=\frac
d{dt}f(x+tv_1)_{|t=0}+i\frac d{dt}f(x+tv_2)_{|t=0} ;$$
on \'ecrira aussi $\partial(v)f$ pour $v\cdot f$.

Dans toute la suite  $G$  d\'esignera un
groupe de Lie r\'eductif connexe muni d'une involution $\sigma$. On note  $H$ la composante connexe de l'\'el\'ement neutre du groupe des
points fixes de $\sigma$,  $\g g$ l'alg\`ebre de Lie de $G$, $\g h$
celle de $H$ et encore $\sigma$ la diff\'erentielle de $\sigma$ en
l'identit\'e. On a alors $\g g=\g h\oplus \g q$ la d\'ecomposition de $\g g$ en
somme de sous-espaces propres associ\'es aux valeurs propres $1$ et $-1$
de $\sigma$, et $H$ agit dans $\g q$ par restriction de la
repr\'esentation adjointe. Les couples $(G,H)$ et $(\g g,\g h)$ sont
appel\'es des paires sym\'etriques r\'eductives.

On utilisera librement les r\'esultats classiques sur les paires
sym\'etriques (voir par exemple \cite{T-Y}).

On fixe une forme bilin\'eaire sym\'etrique non d\'eg\'en\'er\'ee
$\kappa$ sur $\g g$ invariante par $G$ et par $\sigma$. Sa restriction
\`a $\g q$ est alors non d\'eg\'en\'er\'ee. On note de la m\^eme fa\c con son
extension par lin\'earit\'e \`a $\g g_{\C}$.

On notera $\Car(\g q)$ l'ensemble des sous-espaces de Cartan de $\g
q$. Pour $\g a\in \Car(\g q)$, on note $\Delta(\g g_{\C},\g a_{\C})$
l'ensemble des racines de $\g a_{\C}$ dans $\g g_{\C}$. C'est un
syst\`eme de racines non r\'eduit en g\'en\'eral, on note $W(\g a)$, ou
simplement $W$ quand il n'y a pas de risque de confusion, son  groupe de
Weyl. On note $W(H,\g a)$ le quotient du normalisateur de $\g a$  par son
centralisateur dans
$H$ ; c'est un sous-groupe de $W(\g a)$.

L'application de restriction induit un isomorphisme (isomorphisme de Chevalley) de l'alg\`ebre
$\C[\g q_{\C}]^H$ des polyn\^omes $H$-invariants sur $\g q_{\C}$, sur l'alg\`ebre $\C[\g a_{\C}]^{W(\g a)}$ des polyn\^omes $W(\g a)$-invariants sur $\g a_{\C}$.

Soit $\g a\in \Car(\g q)$. On note  $\Delta_r(\g
g_{\C},\g a_{\C})$ le sous-ensemble de $\Delta(\g g_{\C},\g a_{\C})$
form\'ee des racines $\alpha$ telles que $\frac{\alpha}2$ n'est pas une
racine. Le polyn\^ome 
$$\prod_{\alpha\in \Delta_r(\g
g_{\C},\g a_{\C})}\alpha $$
 est invariant par $W(\g a)$, il se prolonge donc de
fa\c con unique en un  polyn\^ome $H$-invariant sur $\g q_{\C}$, qu'on notera
$\it{\Phi}$. Il est clair que $\it{\Phi}$ ne d\'epend pas du choix du
sous-espace de Cartan $\g a$. Si $\g g$ est ab\'elienne, on pose
$\it{\Phi}=1$.

L'alg\`ebre $\C[\g q_{\C}]^H$ est une alg\`ebre de polyn\^omes. On en fixe un
syst\`eme minimal de g\'en\'erateurs homog\`enes $(P_1,\ldots, P_{\ell})$ ,
$\ell$ \'etant le rang de la paire $(\g g,\g h)$, c'est-\`a-dire la
dimension de tout sous-espace de Cartan de $\g q$.

La restriction \`a $\g q_{\C}$ de $\kappa$ est non
d\'eg\'en\'er\'ee ; cela permet de d\'efinir le gradient $\nabla f$ de toute
fonction lisse $f$ sur un ouvert $U$ de $\g q$. Si la fonction est
$H$-invariante, alors $\nabla f$ est un champ de vecteurs
invariant. En particulier on dispose  des champs de vecteurs $\nabla P_i$,
$1\leq i\leq \ell$. Leur restriction \`a tout ouvert de $\g q$ sera
not\'ee de la m\^eme fa\c con. On verra dans la section 8 que la fonction
$\it{\Phi}$
est \`a une constante pr\`es \'egale  \`a 
$$\det (\nabla P_i\cdot P_j)_{1\leq i,j\leq \ell}.$$

Rappelons qu'un ouvert  $\mathcal U$ de $\g q$ est dit {\em
  compl\`etement $H$-invariant} s'il est $H$-invariant et s'il contient
  la composante semi-simple (de la d\'ecomposition de Jordan) de chacun
  de ses \'el\'ements ; on verra dans la section \ref{descente} les
  principales propri\'et\'es de ces ouverts.

Tout champ de vecteurs $X\in \g X(\mathcal U)^H$ induit une d\'erivation
$\C$-lin\'eaire de
$C^{\infty}(\mathcal U)^H$, qu'on notera $D_X$. Le r\'esultat principal
de cet article est le th\'eor\`eme suivant.

\begin{theo}\label{principal}{Soit $\mathcal U\subset \g q$ un ouvert compl\`etement
  $H$-invariant  et soit $X\in \g X(\mathcal U)^H$. Alors
    $D_X$ laisse stable l'id\'eal $\it{\Phi} C^{\infty}(\mathcal U)^H$ de
    $C^{\infty}(\mathcal U)^H$. R\'eciproquement, pour toute
    d\'erivation $\C$-lin\'eaire  $D$
    de $C^{\infty}(\mathcal U)^H$
    qui laisse stable l'id\'eal $\it{\Phi} C^{\infty}(\mathcal U)^H$, il
    existe des fonctions  $\varphi_1\ldots,\varphi_{\ell}\in
    C^{\infty}(\mathcal U)^H$, uniquement d\'etermin\'ees, telles que
    $D=D_X$, o\`u $X=\sum_{1\leq i\leq \ell}\varphi_i\nabla P_i$.}
\end{theo}

La preuve de la premi\`ere partie du th\'eor\`eme sera donn\'ee dans la
section \ref{dem1} et celle de la deuxi\`eme partie dans la section \ref{dem2}.

\begin{rem}{On  peut choisir les $P_i$ \`a valeurs r\'eelles sur $\g q$ ;
    ce fait est bien connu, on peut le d\'emontrer en utilisant par
    exemple les arguments de (\cite{B}, Annexe). Il est facile de voir
    que $\it{\Phi}$ est r\'eelle sur $\g q$. Ceci not\'e,  on peut
    ais\'ement obtenir une
    version r\'eelle du th\'eor\`eme \ref{principal}. Plus pr\'ecis\'ement si
    l'on note $C^\infty_{\R}(\mathcal U)^H$ le sous-espace de
    $C^\infty(\mathcal U)^H$ form\'ee des fonctions r\'eelles, et $\g
    X_{\R}(\mathcal U)^H$ le sous-espace de $\g X(\mathcal U)^H$ form\'e
    des champs de vecteurs \`a valeurs dans $\g q$, alors le th\'eor\`eme
    est encore vrai si l'on remplace dans l'\'enonc\'e $\g X(\mathcal
    U)^H$ par $\g X_{\R}(\mathcal
    U)^H$ et $C^\infty(\mathcal U)^H$ par $C^\infty_{\R}(\mathcal U)^H$ ;
    les $\varphi_i$ sont dans ce cas dans $C^\infty_{\R}(\mathcal U)^H$.}
\end{rem}

\section{D\'eveloppement de Taylor}\label{taylor}

Pour notre formulation des s\'eries de Taylor, il est commode
d'introduire la d\'efinition suivante. Soit $V$ un espace vectoriel
r\'eel de dimension finie. Si $a$ est
un \'el\'ement de $V$, on dit qu'un polyn\^ome $P$ sur $V$ est {\em homog\`ene de
  degr\'e $k$}, $k\in \N$, en $x-a$ si 
$$P(a+t(x-a))=t^kP(x)\quad\mbox{\rm pour tout } x\in V\;\; \mbox{\rm
  et tout } t\in\R ;$$
les polyn\^omes homg\`enes en $x-a$ de degr\'e $0$ sont les
  constantes. L'ensemble de ces polyn\^omes sera not\'e $\C[V_a]_k$.

On note $\C[[V_a]]=\prod_k\C[V_a]_k$ ;
c'est une $\C$-alg\`ebre unitaire  pour le produit :
$$(P_k)_{k\in\N}\cdot(Q_k)_{k\in\N}=(R_k)_{k\in\N}, \mbox{\rm o\`u }
R_k=\sum_{j=0}^kP_jQ_{k-j} ;$$
l'\'el\'ement unit\'e est la suite $(P_k)_{k\in\N}$, avec $P_0=1$ et $P_k=0$
pour $k\geq 1$. Un \'el\'ement $(P_k)_{k\in\N}$ de $\C[[V_a]]$ est
inversible si et seulement si $P_0\neq 0$.

Pour $f\in C^{\infty}(U)$, on appelle {\em s\'erie de Taylor} de $f$ en
$a\in U$, qu'on note $T_af$, l'\'el\'ement  $(P_k)_{k\in\N}$ de $\C[[V_a]]$  d\'efini par
$$P_k(x)=\frac 1{k!}\frac{d^k}{dt^k}f(a+t(x-a))_{|t=0}.$$

 Pour
tout ouvert $U$ contenant $a$, l'application $ T_a$ est un morphisme d'alg\`ebres  de $C^{\infty}(U)$ sur $\C[[V_a]]$, surjectif
d'apr\`es le lemme de Borel. Sa restriction \`a l'alg\`ebre des
 polyn\^omes est un isomorphisme de cette alg\`ebre sur la
sous-alg\`ebre de  $\C[[V_a]]$ form\'ee des suites dont tous les termes
sont nuls sauf un nombre fini. Pour cette raison, il est commode
d'identifier un polyn\^ome avec sa s\'erie de Taylor. 

Le th\'eor\`eme suivant est fondamental pour cet article (voir \cite{Ma},
proposition 4.1).

\begin{theo}\label{malgrange}{ Soient $U$ un ouvert de $V$ et $P$ un polyn\^ome non nul sur $V$. Pour   que $f\in C^{\infty}(U)$ soit divisible par $P$, i.e. il existe    $g\in C^{\infty}(U)$ tel que $f=Pg$,  il faut et il suffit que  pour tout $a\in U$, la s\'erie de Taylor $T_af$ soit divisible par  $P$ dans $\C[[V_a]]$.}
\end{theo}

L'action de $V_{\C}$ dans  $C^\infty(V)$ par
d\'erivations  se prolonge
de fa\c con unique en une action de  $S(V_{\C})$, qui identifie cette
alg\`ebre \`a   l'alg\`ebre des op\'erateurs diff\'erentiels \`a coefficients
constants sur $V$ ; on note $\partial(u)$ l'op\'erateur diff\'erentiel
associ\'e \`a $u\in S(V_{\C})$. Pour  $f\in C^{\infty}(U)$,  $a\in U$ et  $u\in S(V_{\C})$, on a
$$\partial(u)f(a)=\sum_k\partial (u)P_k(a) ;$$
o\`u $T_af=(P_k)_{k\in \N}$.

Soit $a\in U$. Tout vecteur $v\in V_{\C}$ d\'efinit une d\'erivation $\widehat{\partial}(v)$ de
$\C[[V_a]]$ par
$$\widehat{\partial}(v)\cdot((P_k)_{k\in\N})=(Q_k)_{k\in \N}, \quad \mbox{\rm avec }
Q_k=\partial(v)P_{k+1}.$$

Il est clair que pour tout $f\in C^{\infty}(U)$ et pour tout $v\in V_{\C}$, on a
$$T_a(\partial(v)\cdot f)=\widehat{\partial}(v)\cdot (T_af).$$

Si $X\in \g X(U)$, il existe
une unique d\'erivation de $\C[[V_a]]$, qu'on note $T_aX$, d\'efinie de
la fa\c con suivante. Si on \'ecrit  $X=\sum_{i}{\varphi_iv_i}$, $\varphi_i\in C^{\infty}(U)$, alors
$$T_aX=\sum (T_a\varphi_i) \widehat{\partial}(v_i).$$
On a
$$(T_aX)(T_af)=T_a(Xf),\quad\mbox{\rm pour tout } f\in C^{\infty}(U).$$

L'application $\sum_{i}S_i\otimes v_i\mapsto \sum_iS_i\widehat{\partial}(v_i)$ est un
isomorphisme lin\'eaire de $\C[[V_a]]\otimes_{\C}V_{\C}$ sur l'espace
vectoriel des d\'erivations $\C$-lin\'eaires de $\C[[V_a]]$ (voir
\cite{Bour}, §4, no. 6). Avec cette
identification, si $X\in  \g X(U)$, alors $T_aX$ s'identifie \`a la
famille $(D_k)_{k\in\N}$ de polyn\^omes homog\`enes sur $V$ \`a valeurs dans
$V_{\C}$ d\'efinie par
$$D_k(x)=\frac 1{k!}\frac{d^k}{dt^k}X(a+t(x-a))_{|t=0}.$$

Si $V_{\C}$ est muni d'une forme bilin\'eaire sym\'etrique non d\'eg\'en\'er\'ee, ce
qui permet de d\'efinir le gradient $\nabla f$ d'une fonction lisse $f$
sur $U$, et si $R$ est un polyn\^ome homog\`ene en $x-a$ de degr\'e $d\geq
1$, on a 
$$(T_a \nabla R)\cdot ((P_k)_{k\in \N})=(Q_k)_{k\in \N}, \quad \mbox{\rm avec }
  Q_k=\left\{\begin{array}{ll}\nabla R\cdot P_{k-d+2}& \mbox{\rm si  } k\geq d-1\\  0&\mbox{\rm si  } k< d-1\end{array}\right..$$

Si $f$ est une fonction sur un ensemble $A$ et si $B$ est une partie
de $A$, on note $f|B$ la restriction de $f$ \`a $B$.

Soit $F$ un sous-espace vectoriel de $V$ et soit $a\in U\cap F$. Pour
tout $f\in C^{\infty}(U)$,  si l'on note $T_af=(P_k)_{k\in \N}$, alors
il est clair d'apr\`es les d\'efinitions que
$$T_a(f|(U\cap F))=(P_k|F)_{k\in \N}.$$

\section{M\'ethode de descente}\label{descente}

Nous rappelons dans cette section les principaux outils de la m\'ethode
de descente de Harish-Chandra. 

Si $a\in\g q$ est semi-simple, on note $G^a$ (resp. $H^a$) le fixateur
de $a$ dans $G$ (resp. $H$), le groupe $H^a$ n'est pas forc\'ement connexe. On note $\g g^a$,  $\g h^a$ et $\g q^a$  les centralisateurs de
$a$ respectivement  
dans $\g g$,  $\g h$ et $\g q$. Alors $(\g g^a,\g h^a)$ est une paire
sym\'etrique r\'eductive et $\g g^a=\g h^a\oplus \g q^a$. Un ouvert
$\mathcal V$ de $\g q^a$ st dit $H^a$-{\it compl\`etement invariant} s'il est
$H^a$-invariant et s'il contient la composante semi-simple de chacun
de ses \'el\'ements.

\newcounter{bean}
Soit $\mathcal U$ un ouvert  compl\`etement $H$-invariant de $\g
q$ et soit $a\in \mathcal U$ semi-simple. On a (\cite{Ka},
proposition 3.1) :

\begin{prop}\label{desprop}{Il existe un voisinage ouvert $\mathcal V$ de $a$
dans $\g q^a$ compl\`etement $H^a$-invariant et 
v\'erifiant :
\begin{list}
{\mbox{\rm (\roman{bean})}}{\usecounter{bean}
\setlength{\rightmargin}{\leftmargin}} 
\item {l'application $\pi : H\times V\rightarrow \g q$ d\'efinie
  par  $\pi(h,y)=h\cdot y$ est une submersion et $\mathcal
  W=\pi(H\times V)$ est un ouvert compl\`etement $H$-invariant de $\g q$
  inclus dans $\mathcal U$ ;}

\item { l'application de restriction est un isomorphisme de
  $C^{\infty}(\mathcal W)^H$ sur $C^{\infty}(\mathcal V)^{H^a}$}.
\end{list}}
\end{prop}

Il est \`a noter que l'application $\pi$ induit une
bijection de l'espace fibr\'e $H\times_{H^a}\mathcal V$ sur l'ouvert
$\mathcal W$. Donc si $(\rho, E)$ est une repr\'esentation de dimension
finie de $H$, toute
fonction $f\in C^{\infty}(\mathcal V,E)^{H^a}$ se prologe de fa\c con
unique en une fonction $\widetilde{f}\in C^{\infty}(\mathcal W,E)^{H}$.

\begin{cor}\label{descor}
{Il existe une fonction $\chi\in
C^{\infty}(\g q)^H$ telle que $0\leq \chi(y)\leq 1$ pour tout $y\in \g
q$, $\chi=1$ dans un voisinage de $a$ et le support de $\chi$ est inclus
dans $\mathcal U$.}
\end{cor}

 La preuve est une adaptation au cas sym\'etrique de
celle du th\'eor\`eme 1.27 de (\cite{Va}, Part I) ; tous les ingr\'edients
se trouvent dans la preuve de (\cite{Ka},
proposition 3.1).

\begin{lemme}\label{restcv}{Soit $f\in C^{\infty}(\mathcal U)^H$ et soit
    $a\in\mathcal U$ semi-simple. Alors $\nabla f(y)\in \g q_{\C}^a$ pour
    tout $y\in \mathcal U\cap \g q^a$.}
\end{lemme}
\begin{dem} La restriction de $\kappa$ \`a $\g q_{\C}^a$ 
 \'etant  non d\'eg\'en\'er\'ee, il suffit de montrer que la diff\'erentielle
 $df(y)$ de $f$ en $y$ est nulle sur l'orthogonal $\g m$ de $\g q^a$
 dans $\g q$ pour $\kappa$. 

On note ${'}\g q^a$ l'ensemble des $y\in\g q^a$ tels que $\det (\ad y)_{\g
  g/\g g^a}\neq 0$ ; c'est un ouvert dense de $\g q^a$. Pour tout
  $y\in\ {'}\g q^a$, on a $\g m\subset [y,\g h]$.

Soit $\xi\in \g h$ et soit $y\in  \mathcal U\cap\ {'}\g q^a$. La fonction
$t\mapsto f(\exp t\xi\cdot h)$ est constante sur $\R$, donc
$$df(y)\cdot [y,\xi]=\frac d{dt}f((\exp -t\xi\cdot h)_{|t=0}=0.$$
Il s'ensuit que $df(y)$ est nulle sur $\g m$ pour tout $y\in \mathcal
U\cap\ {'}\g q^a$, et donc, par  densit\'e de $\mathcal U\cap\ {'}\g q^a$ dans $\mathcal
U\cap \g q^a$, $df(y)$ est nulle sur $\g m$ pour tout $y\in \mathcal
U\cap \g q^a$.
\end{dem}

Si  $\mathcal A$ est une $\C$-alg\`ebre associative et unitaire, on note $\Der(\mathcal
A)$ l'ensemble de ses d\'erivations $\C$-lin\'eaires.

\begin{prop}\label{locder}{Soit $\mathcal U'\subset \mathcal U$ deux
    ouverts compl\`etement  $H$-invariants de $\g q$ et soit $D\in\Der
    (C^{\infty}(\mathcal  U)^H)$. Alors il existe une unique d\'erivation 
    $D'\in\Der( C^{\infty}(\mathcal U')^H)$ telle que
$$D'\cdot (f|\mathcal U')=(D\cdot f)|\mathcal U',\quad \mbox{\rm pour
      tout } f\in C^{\infty}(\mathcal U)^H.$$}
\end{prop}
\begin{dem} Tenant compte du corollaire \ref{descor}, la preuve est
      exactement la m\^eme que dans le cas usuel (sans action de
      groupe), voir par exemple (\cite{La}, th\'eor\`eme III.17).

\end{dem}

\section{Division par un polyn\^ome invariant}\label{division}

L'involution $\sigma$ de $\g g$ d\'efinit par lin\'earit\'e une involution
$\sigma_{\C}$ de $\g g_{\C}$.  On note $H_{\C}$ la composante neutre
du commutant $\sigma_{\C}$ dans le  groupe adjoint de
$\g g_{\C}$. Ce groupe op\`ere
naturellement dans $\g q_{\C}$, et les sous-espaces de Cartan de $\g
q_{\C}$ sont tous conjugu\'es sous son action.

Soit $\mathcal U$ un ouvert
    compl\`etement $H$-invariant de $\g q$. Pour all\'eger les notations
    dans ce qui suit, pour $f\in
    C^{\infty}(\mathcal U)^H$ et pour $\g a\in \Car(\g q)$, on note
    $f|\g a$ pour $f|(\mathcal U\cap \g a)$.

\begin{prop}\label{divisioninvariante}{Soient  $f\in C^{\infty}(\mathcal U)^H$ et $P\in \C[\g q]^H$ non
    nul. On suppose que pour tout $\g a\in \Car(\g q)$, il existe
    $g_{\g a}\in C^{\infty}(\mathcal U\cap\g a)$ tel que $f{|\g a}=(P{|\g a})g_{\g
    a}$. Alors il existe $g\in C^{\infty}(\mathcal U)^H$ tel que $f=Pg$.}
\end{prop}

\begin{dem} On commence par montrer qu'il existe une fonction continue
    $g$ sur $\mathcal U$, forc\'ement unique, telle que $g{|\g a}=g_{\g a}$
    pour tout $\g a\in\Car(\g q)$.  

Il est clair que pour tout $h\in
    H$ et pour tout $\g a\in \Car(\g q)$, on a $h\cdot g_{\g
    a}=g_{h\cdot \g a}$. Donc, d'apr\`es (\cite{Ka}, lemme 6.1), il
    suffit de montrer que pour tous $\g a, \g b\in\Car(\g q)$ et pour
    tout $a\in \mathcal U\cap \g a\cap\g b$, on a $g_{\g a}(a)=g_{\g b}(a)$.

Supposons donc $\g a,\g b, a$ donn\'es comme ci-dessus. Il existe $u\in
S(\g a_{\C})$ homog\`ene tel que $\partial(u)(P{|\g a})(a)=1$ et
$\partial(v)(P{|\g a})(a)=0$ pour tout $v\in S(\g a_{\C})$ homog\`ene de
degr\'e strictement plus petit que celui de $u$ (le degr\'e de $u$ est le
degr\'e du premier terme non nul dans le d\'eveloppement de Taylor de
$P{|\g a}$
en $a$). Alors
\begin{equation}\label{un}
\partial(u)(f{|\g a})(a)=\partial(u)((P{|\g a})g_{\g a})(a)=g_{\g
  a}(a).
\end{equation}

On note $H_{\C}^a$ le sous-groupe des \'el\'ements de $H_{\C}$ qui fixent
$a$. Il existe $h\in H_{\C}^a$ tel que $h\cdot \g a_{\C}=\g
b_{\C}$. D'apr\`es (\cite{Ka}, th\'eor\`eme 5.1), 
\begin{equation}\label{deux}
\partial(h\cdot u)(f{|\g b})(a)=\partial(u)(f{|\g a})(a)\quad\mbox{\rm
  et} \quad\partial(h\cdot u)(P{|\g b})(a)=\partial(u)(P{|\g a})(a)=1 ;
\end{equation}
de plus, pour tout $v\in S(\g b_{\C})$ homog\`ene de degr\'e strictement plus petit
que celui de $u$, on a 
$$\partial(v)(P{|\g b})(a)=\partial(h^{-1}\cdot v)(P{|\g a})(a)=0,$$
donc  
\begin{equation}\label{trois}
\partial(h\cdot u)(f{|\g b})(a)=\partial(h\cdot u)((P{|\g b})g_{\g b})(a)=g_{\g
  b}(a).
\end{equation}
Il d\'ecoule alors de (\ref{un}),(\ref{deux}) et (\ref{trois})  que $g_{\g a}(a)=g_{\g b}(a)$.

Pour montrer que $g$ est lisse, on utilise (\cite{Ka}, th\'eor\`eme
5.1). La fonction $g$ doit v\'erifier deux conditions. La premi\`ere est
l'invariance de $g_{\g a}$ par le groupe $W(H,\g a)$ pour tout $\g a\in \Car(\g
q)$ ; cette propri\'et\'e d\'ecoule facilement de l'invariance de $f{|\g
  a}$ et de $P{|\g a}$ par ce groupe. La seconde est que pour tout
$a\in\g q$ semi-simple, pour tous $\g a_1, \g a_2\in \Car(\g q)$
contenant $a$,  pour tout $h\in H_{\C}^a$ tel que $h\cdot \g
{a_1}_{\C}= \g {a_2}_{\C}$ et pour tout $u\in S(\g {a_1}_{\C})$, on a
\begin{equation}\label{recollement}
\partial(u)g_{\g a_1}(a)=\partial(h\cdot u)g_{\g a_2}(a).
\end{equation}

Supposons que nous ayons toutes ces donn\'ees. On fixe une base
$(v_1,\ldots,v_{\ell})$ de $\g a_1$, et pour
$\alpha=(\alpha_1,\ldots,\alpha_{\ell})\in\N^{\ell}$, on note
$v^{\alpha}=v_1^{\alpha_1}\cdots v_{\ell}^{\alpha_{\ell}}$, produit
dans $S(\g {a_1}_{\C})$. La famille $(v^{\alpha})_{\alpha\in\N^{\ell}}$
est une base de $S(\g {a_1}_{\C})$. On va donc montrer
(\ref{recollement}) pour les $v^{\alpha}$. Pour cela nous avons besoin
d'une relation d'ordre sur $\N^{\ell}$ pour faire une r\'ecurrence.

Si $\alpha=(\alpha_1,\ldots,\alpha_{\ell})\in \N^{\ell}$, on note
$\vert \alpha\vert=\alpha_1+\cdots+\alpha_{\ell}$ sa longueur. On
dispose de la relation d'ordre sur $\N^{\ell}$ :
$\alpha\leq\beta$ si $\alpha_1\leq
\beta_1,\ldots,\alpha_{\ell}\leq\beta_{\ell}$, de la relation d'ordre
lexicographique : $\alpha\leq_L\beta$ si $\alpha=\beta$ ou s'il existe
$1\leq j\leq \ell$ tel que $\alpha_i=\beta_i$ pour $1\leq i<j$ et
$\alpha_j<\beta_j$. 
 On d\'efinit une nouvelle
relation d'ordre total sur $\N^{\ell}$ par
$$\alpha\preceq\beta \quad\mbox{ \rm si } \quad
\vert\alpha\vert<\vert\beta\vert \mbox{ \rm ou }
\vert\alpha\vert=\vert\beta\vert  \mbox{ \rm et  } \alpha\leq_L\beta.$$ 
Si $\alpha\preceq\beta$ et $\alpha\neq \beta$, on \'ecrit $\alpha \prec\beta$.

On va donc montrer
(\ref{recollement}) pour les $v^{\alpha}$ par r\'ecurrence suivant l'ordre $\preceq$. On l'a d\'ej\`a  prouv\'e plus haut pour $\vert \alpha\vert=0$. 

On note $d$ le plus petit entier $k$ tel que la restriction \`a $\g a_1$
de la composante homog\`ene en
$x-a$ de degr\'e $k$ de $T_aP$ soit non nulle. Soit $\alpha^0$ le plus
petit $\alpha\in\N^{\ell}$ pour l'ordre $\preceq$ tel que
$\partial(v^{\alpha})(P{|\g a_1})(a)\neq 0$. On a
$\vert\alpha^0\vert=d$, et, pour tout $\gamma \in\N^{\ell}$ de longueur
$<d$, $\partial(v^{\gamma})(P{|\g a_1})(a)=0$. 

D'apr\`es (\cite{Ka}, th\'eor\`eme 5.1)), pour toute fonction $F\in
C^{\infty}(\mathcal U)^H$ et pour tout $u\in S({\g a_1}_{\C}$, on a
\begin{equation}\label{recollementbis}
\partial(u)(F|\g a_1)(a)=\partial(h\cdot u)(F|\g
  a_2)(a).
\end{equation}
Il d\'ecoule donc de ce qui pr\'ec\`ede que  $\partial(h\cdot
v^{\alpha^0})(P{|\g a_2})(a)\neq 0$ et $\partial(h\cdot
v^{\gamma})(P{|\g a_2})(a)= 0$ si $\vert\gamma\vert<d$.

Soit $\beta\in\N^{\ell}$ tel que $\vert\beta\vert\geq 1$. On suppose
que la condition (\ref{recollement}) est v\'erifi\'ee pour les
$v^{\gamma}$, $\gamma\prec\beta$. Avec les notations
usuelles, on a
\begin{eqnarray}
\partial(v^{\beta+\alpha^0})f_{\g a_1}(a)&=&\sum_{\gamma\in\N^{\ell}\atop
  \gamma\leq \beta+\alpha^0}\binom{\beta}{\gamma}\partial(v^{\gamma})(P{|\g
  a_1})(a)\partial(v^{\beta+\alpha^0-\gamma})g_{\g a_1}(a)\label{eq1}\\
\partial(h\cdot v^{\beta+\alpha^0})f_{\g
  a_2}(a)&=&\sum_{\gamma\in\N^{\ell}\atop\gamma\leq
  \beta+\alpha^0}\binom{\beta}{\gamma}\partial(h\cdot
  v^{\gamma})(P{|\g  a_2})(a)\partial(h\cdot v^{\beta+\alpha^0-\gamma})g_{\g a_2}(a)\label{eq2}
\end{eqnarray}

Il d\'ecoule de (\ref{recollementbis}), appliqu\'e \`a $f$,  que
les membres de gauche de (\ref{eq1}) et de (\ref{eq2}) sont
\'egaux. Pour les membres de droite, on a  $\partial(v^{\gamma})(P{|\g
  a_1})(a)=\partial(h\cdot v^{\gamma})(P{|\g
  a_2})(a)$ pour tout $\gamma$ ;  ces termes sont nuls si $\vert
\gamma\vert<d$ et l'hypoth\`ese de r\'ecurrence implique que $\partial(v^{\beta+\alpha^0-\gamma})g_{\g
  a_1}(a)=\partial(h\cdot v^{\beta+\alpha^0-\gamma})g_{\g a_2}(a)$ si
$\beta+\alpha^0-\gamma\prec\beta$ ; cela est vrai en particulier pour
tout  $\gamma$ de longueur $\vert \gamma\vert>d$. Il s'ensuit que les
sommes des termes correspondant aux $\gamma$ de longueur $d$ dans les
membres de droite de (\ref{eq1}) et de (\ref{eq2}) sont \'egales, soit
$$\sum_{
  \gamma\leq \beta+\alpha^0; \vert\gamma\vert=d}\binom{\beta}{\gamma}\partial(v^{\gamma})(P{|\g
  a_1})(a)\left[\partial(v^{\beta+\alpha^0-\gamma})g_{\g a_1}(a)-\partial(h\cdot v^{\beta+\alpha^0-\gamma})g_{\g a_2}(a)\right]=0$$

Dans
cette somme la contribution des $\gamma\prec \alpha^0$ est nulle,
puisque $\partial(v^{\gamma})(P{|\g a_1})(a)=0$ par d\'efinition de
$\alpha^0$. Il en est de m\^eme de la contribution des $\gamma$
v\'erifiant $\beta+\alpha^0-\gamma\prec \beta$, d'apr\`es l'hypoth\`ese de
r\'ecurrence. Il ne reste plus que la contribution des $\gamma$
v\'erifiant :
$$\gamma\leq \beta+\alpha^0, \quad\vert\gamma\vert=d,\quad\alpha^0\preceq \gamma,\quad
\beta\preceq \beta+\alpha^0-\gamma. 
$$
Un tel $\gamma$ v\'erifie $\vert\gamma\vert=\vert\alpha^0\vert$, et donc
$\vert\beta+\alpha^0-\gamma\vert=\vert\beta\vert$ ; d'o\`u
$\beta\leq_L\beta+\alpha^0-\gamma$ et $\alpha^0\leq_L\gamma$. Seul
$\gamma=\alpha^0$ est solution de ces deux in\'egalit\'es. D'o\`u
$$ \partial(v^{\alpha^0})(P{|\g
  a_1})(a)\left[\partial(v^{\beta})g_{\g a_1}(a)-\partial(h\cdot
  v^{\beta})g_{\g a_2}(a)\right]=0.$$
Comme $ \partial(v^{\alpha^0})(P{|\g a_1})(a)\neq 0$, il s'ensuit que 
$$\partial(v^{\beta})g_{\g a_1}(a)=\partial(h\cdot
  v^{\beta})g_{\g a_2}(a).$$
\end{dem}

\section{D\'erivations dans les sous-espaces de Cartan}\label{localisationder}

On fixe un sous-espace de Cartan $\g a\in\Car(\g q)$ et un \'el\'ement
$a\in \g a$. On reprend les notations de la section
 \ref{descente} et on note $\g s$ (resp. $\g c$) l'intersection de l'alg\`ebre
 d\'eriv\'ee (resp. du centre) de $\g g^a$ avec $\g q^a$  de sorte que $\g
 q^a=\g s\oplus \g c$. On note $\g b=\g a\cap \g s$, alors $\g a=\g
 b\oplus \g c$.

 On note $W^a$ l'ensemble des \'el\'ements de
$W(\g a)$ qui fixent $a$. C'est un sous-groupe de $W(\g a)$ qui s'identifie au groupe de Weyl du syst\`eme de
 racines $\Delta(\g g_{\C}^a,\g a_{\C})$. Il op\`ere trivialement sur
 $\g c$.

L'alg\`ebre des invariants  $\C[\g a]^{W^a}$ est \'egale \`a $\C[\g
  c]\otimes \C[\g b]^{W^a}$, et $ \C[\g b]^{W^a}$ est une
  alg\`ebre de polyn\^omes. On  fixe un syst\`eme minimal de
  g\'en\'erateurs homog\`enes $q_1\ldots,q_r$, $r=\dim\g b$, de $\C[\g
  b]^{W^a}$, qu'on regarde
  comme des polyn\^omes sur $\g a$ invariants par translation par les
  \'el\'ements de $\g c$. Sans nuire \`a la g\'en\'eralit\'e, on suppose que
  les $q_i$ sont \`a valeurs r\'eelles sur $\g a$. On fixe des
  coordonn\'es lin\'eaires r\'eelles $x_1,\ldots, x_{\ell-r}$
sur $\g c$, qu'on regarde comme des formes lin\'eaires sur $\g a$ nulles
  sur $\g b$. Alors les polyn\^omes $q_1,\ldots,q_{\ell}$, avec
  $q_{r+j}=x_j-a_j$ pour $1\leq j\leq \ell-r$, sont homog\`enes en $x-a$
  et engendrent $\C[\g a]^{W^a}$. On notera $d_i$ le degr\'e
  d'homog\'en\'eit\'e de $q_i$.

On note $\C[[q_1,\ldots,q_{\ell}]]$ l'alg\`ebre des s\'eries formelles en
les variables $q_1,\ldots,q_{\ell}$. Si l'on note, pour
$\alpha=(\alpha_1,\ldots,\alpha_{\ell})\in\N^{\ell}$,
$q^{\alpha}=q_1^{\alpha_1}\cdots q^{\alpha_{\ell}}$, alors tout
\'el\'ement de $\C[[q_1,\ldots,q_{\ell}]]$ s'\'ecrit de fa\c con unique sous la forme
$$\sum_{\alpha\in\N^{\ell}}c_{\alpha}q^{\alpha},
\quad c_{\alpha}\in\C.$$

Soit $\mathcal U$  un ouvert compl\`etement $H$-invariant de $\g q$
contenant $a$. On note $\mathcal V=\mathcal U\cap \g a$ et
$C^{\infty}(\mathcal V)^{\inv}$
l'ensemble des restrictions  $f{|\mathcal V}$, $f\in C^{\infty}(\mathcal U)^H$ ;
on trouvera une caract\'erisation de cet ensemble dans
(\cite{Ka2}). En particulier, pour tout $f\in C^{\infty}(\mathcal
V)^{\inv}$,  la s\'erie de taylor de $f$ en $a$
appartient \`a l'alg\`ebre $\C[[\g a_a]]^{W^a}$  des invariants par $W^a$ dans
$\C[[\g a_a]]$.

L'application qui associe \`a un \'el\'ement
$\sum_{\alpha\in\N^{\ell}}c_{\alpha}q^{\alpha}$ de  $\C[[q_1,\ldots,q_{\ell}]]$
la suite
$(P_k)_{k\in\N}\in\C[[\g a_a]]$ d\'efinie par
$$P_k=\sum_{\alpha\in \N^{\ell}\atop
  \alpha_1d_1+\cdots\alpha_{\ell}d_{\ell}=k}c_{\alpha} q^{\alpha}$$
est un isomorphisme d'alg\`ebres de $\C[[q_1,\ldots,q_{\ell}]]$ sur
  $\C[[\g a_a]]^{W^a}$. Dans la suite on identifiera ces deux alg\`ebres.

\begin{lemme}\label{taylorinv}{L'application
  $f\mapsto T_af$ de
  $C^{\infty}(\mathcal V)^{\inv}$ dans $\C[[\g a_a]]^{W^a}$  est un morphisme de $\C$-alg\`ebres surjectif.
}
\end{lemme}
\begin{dem} Soit $\eta$ l'application de $\mathcal V$ dans  $\R^{\ell}$
  d\'efinie par $\eta(x)=(q_1(x),\ldots,q_{\ell}(x))$. On a $\eta(a)=0$. Soit
  $S=\sum_{\alpha\in\N^{\ell}}c_{\alpha}q^{\alpha}\in
  \C[[q_1,\ldots,q_{\ell}]]$. On note $y_1,\ldots,y_{\ell}$ les
  coordonn\'ees de $\R^{\ell}$. D'apr\`es le lemme de Borel, il existe
  $\varphi\in C^{\infty}(\R^{\ell})$ tel que 
$$\frac{\partial^{\alpha_1}}{\partial y_1^{\alpha_1}}\cdots \frac{\partial^{\alpha_{\ell}}}{\partial y_{\ell}^{\alpha_{\ell}}}\varphi(0)=\alpha!
  c_{\alpha},\quad\mbox{\rm pour tout }\alpha\in\N^{\ell}.$$ 
 Un calcul simple montre que $T_a(\varphi\rond \eta)=S$.

 Le groupe $H_{\C}^a$ n'est pas connexe en g\'en\'eral, cependant il
  d\'ecoule de (\cite{K-R}, proposition 10 et la discussion page 780) que les
  polyn\^omes invariants par $H_{\C}^a$ ou par sa composante neutre
  sont les m\^emes, d'o\`u
  l'isomorphisme de Chevalley $\C[\g q^a]^{H_{\C}^a}\simeq \C[\g
  a]^{W^a}$, et comme l'action de $H^a$ dans $\g q^a$ se factorise \`a travers
  $H_{\C}^a$, on a $\C[\g q^a]^{H_{\C}^a}=\C[\g q^a]^{H^a}$. Donc,
  pour tout $1\leq i\leq \ell$, il existe un unique \'el\'ement $Q_i$ de
 $\C[\g q^a]^{H^a}$ dont la restriction \`a $\g a$ est \'egal \`a $q_i$ ; ce
  polyn\^ome est \`a valeurs r\'eelles sur $\g q^a$.
  D'apr\`es la proposition \ref{desprop}, les $Q_i$ se prolongent en des
  fonctions $H$-invariantes $\widetilde{Q_i}$ sur un voisinage ouvert 
  compl\`etement $H$-invariant $\mathcal U'$ de $a$. La fonction
  $f(x)=\varphi(\widetilde{Q_1}(x),\ldots,\widetilde{Q}_{\ell}(x))$ appartient \`a $C^{\infty}(\mathcal
  U')^H$, et quitte \`a la multiplier par une fonction plateau
  (corollaire \ref{descor}), on peut la supposer dans $C^{\infty}(\mathcal
  U)^H$. Sa restriction \`a $\mathcal V$ co\" incide avec $\varphi\rond \eta$ dans
  un voisinage de $a$, donc $T_a(f|\mathcal V)=S$.
\end{dem}

\begin{lemme}\label{derpoint}{Soit $\delta\in\Der(C^{\infty}(\mathcal
  V)^{\inv})$. Il existe une unique d\'erivation $\delta_a\in\Der(\C[[\g  a_a]]^{W^a})$ telle que
$$T_a(\delta f)=\delta_a(T_af),\quad \mbox{\rm pour tout } f\in
  C^{\infty}(\mathcal  V)^{\inv}.$$}
\end{lemme}
\begin{dem} Il suffit de montrer que si la s\'erie de Taylor en $a$ de
  $f\in C^{\infty}(\mathcal V)^{\inv}$  est nulle, alors il en est de
  m\^eme de $T_a(\delta f)$. 

Soit donc $f\in  C^{\infty}(\mathcal V)^{\inv}$ tel que $T_af=0$. Si $\chi\in C^{\infty}(\mathcal U)^H$ est \'egale \`a
$1$ dans un voisinage de $a$, alors 
$$T_a(\delta
(f(\chi|_{\mathcal V})))= T_a(\delta f)T_a(\chi|_{\mathcal V})+T_a(f)T_a(\delta( \chi|_{\mathcal V}))= T_a(\delta f),$$
car $T_af=0$ et  $T_a(\chi|_{\mathcal V})$ est l'\'el\'ement unit\'e de $\C[[\g
    a_a]]^{W^a}$. Cela nous permettra au besoin de remplacer $f$ par une
fonction de la forme $f(\chi|_{\mathcal V})$ sans nuire \`a la g\'en\'eralit\'e.

Pour $x\in\g a_{\C}$, on pose $q(x)=q_1(x)^2+\cdots
  +q_{\ell}(x)^2$. Alors $q$ est un polyn\^ome  $W^a$-invariant \`a valeurs positives sur  $\g a$,
  et $a$ est l'unique \'el\'ement de $\g a$ o\`u il s'annule. Avec les
  notations de la preuve du lemme \ref{taylorinv}, on
  note $Q=Q_1^2+\cdots+Q_{\ell}^2$. Alors $Q$ est \`a valeurs positives
  sur $\g q^a$ et le seul point o\`u il s'annule est $a$. On note
  $\widetilde{Q}$ son prolongement $H$-invariant \`a $\mathcal U'$, alors
  les seuls points de $\mathcal U'$ o\`u $\widetilde{Q}$ s'annule sont les
  \'el\'ements de l'orbite de $a$ sous $H$. 

 Soit $\widetilde{f}$ un
  prolongement de $f$ en un \'el\'ement de $C^{\infty}(\mathcal
  U)^H$. L'hypoth\`ese faite sur $f$ implique (\cite{Ka}, th\'eor\`eme 5.1)
  que pour tout $\g b\in\Car(\g q)$ et tout $x\in \g b\cap H\cdot a$,
  la s\'erie de Taylor de $\widetilde{f}|(\mathcal U\cap \g b)$ en $x$ est
  nulle. Quitte \`a multiplier $\widetilde{f}$ par
  une fonction $\chi$ comme ci-dessus, on peut supposer que
  $\widetilde{f}$ est \`a support dans $\mathcal U'$. Donc, pour tout $\g
  b\in\Car(\g q)$, la s\'erie de Taylor de $ \widetilde{f}|(\mathcal U\cap
  \g b)$ est nulle en tout point de $\mathcal U\cap
  \g b$ o\`u $\widetilde{Q}|\g b$ s'annule. Il d\'ecoule alors du th\'eor\`eme
  \ref{malgrange} que, pour tout entier $k\geq 1$, le polyn\^ome $\widetilde{Q}^k|\g b$ divise $ \widetilde{f}|(\mathcal U\cap
  \g b)$, et donc, d'apr\`es la proposition  \ref{divisioninvariante},
  $\widetilde{Q}^k$ divise $\widetilde{f}$. Donc,  pour tout entier $k\geq
  1$, il existe $g_k\in C^{\infty}(\mathcal V)^{\inv}$ tel que
  $f=q^kg_k$. Un calcul simple montre que, pour tout $k\in \N$,  
$$\frac{d^k}{dt^k}(\delta (q^{k+2}g_{k+2}))(a+t(x-a))_{|t=0}=0.$$
il en d\'ecoule que $T_a(\delta f)=0$.

\end{dem}

On note $\it{\Phi}^a$ l'analogue de $\it{\Phi}$ pour la paire
sym\'etrique r\'eductive $(\g g^a,\g h^a)$. Il existe un polyn\^ome $H^a$-invariant $\it{\Psi}^a$ sur $\g q^a$ tel que 
$\it{\Phi}|\g q^a=\it{\Psi}^a\it{\Phi}^a$. La restriction de $\it{\Psi}^a$ \`a $\g a$ est
\'egale \`a
$$\prod_{\alpha\in\Delta_r(\g g_{\C},\g a_{\C})\atop \alpha(a)\neq
 0}\alpha, $$
donc $\it{\Psi}^a(a)\neq 0$.   La restriction de $\it{\Phi}^a$ \`a $\g a$ est un
polyn\^ome homog\`ene en $x-a$ invariant par $W^a$ ; on l'identifiera,
ainsi que tout polyn\^ome $W^a$-invariant, \`a son image dans $\C[[\g a_a]]^{W^a}$.

\begin{lemme}\label{idealpreserv}{Soit $\delta\in
\Der(C^{\infty}(\mathcal V)^{\inv})$. On suppose que $\delta$ laisse
stable l'id\'eal $(\it{\Phi}|\g a)C^{\infty}(\mathcal V)^{\inv}$. Alors
$\delta_a$ laisse stable l'id\'eal $(\it{\Phi}^a|\g a)\C[[\g a_a]]^{W^a}$ de
$\C[[\g a_a]]^{W^a}$.}
\end{lemme}
\begin{dem} Compte tenu du lemme \ref{taylorinv}, il est clair que $\delta_a$ laisse stable l'id\'eal de
$\C[[\g a_a]]^{W^a}$ engendr\'e par $\it{\Phi}|\g
  a$. Cet id\'eal est aussi l'id\'eal engendr\'e par $\it{\Phi}^a|\g a$, car
  $\it{\Psi}^a|\g a$ est inversible dans $\C[[\g a_a]]^{W^a}$. 
\end{dem}

 Les champs de vecteurs $\nabla q_i$ induisent des d\'erivations de
$\C[\g a]^{W^a}$ et de
$\C[[\g a_a]]^{W^a}$, qu'on continue de noter $\nabla q_i$. On note
 $\g X_p(\g a)$ l'ensemble
des champs de vecteurs sur $\g a$ \`a coefficients polynomiaux.

\begin{prop}\label{caspolynomial}{\mbox{\rm (i)} Tout $X\in\g X_p(\g a)^{W^a}$ s'\'ecrit de fa\c con
unique $X=\sum_{1\leq i\leq \ell}R_i \nabla q_i$, o\`u les $R_i$ sont
des \'el\'ements de $\C[\g a]^{W^a}$.

\mbox{\rm (ii)} L'application $X\mapsto D_X$ est une bijection de $\g X_p(\g a)^{W^a}$ sur le sous-espace de  $\Der(\C[\g a]^{W^a})$ form\'e des
  d\'erivations qui laissent stable l'id\'eal $(\it{\Phi}^a|\g a)\C[\g a]^{W^a}$.}
\end{prop}

L'assertion (i) d\'ecoule du th\'eor\`eme principal de
\cite{So}. L'assertion (ii) est la formulation de Ra\"is (\cite{Ra}, §
3.2.) d'un cas particulier du r\'esultat
principal de \cite{Sch}.

\begin{prop}\label{derformel}{Toute d\'erivation $\delta$ de $\C[[\g a_a]]^{W^a}$
qui laisse stable l'id\'eal engendr\'e par $\it{\Phi}^a|\g a$ est de la forme
$\sum_{1\leq i\leq \ell}S_i\nabla q_i$, o\`u les $S_i$ sont des \'el\'ements
de $\C[[\g a_a]]^{W^a}$ uniquement d\'etermin\'es.}
\end{prop}
\begin{dem} La preuve est inspir\'ee de celle de la proposition 6.1 de
  \cite{Sch}. On utilise l'identification $\C[[\g a_a]]^{W^a}\simeq
  \C[[q_1,\ldots,q_{\ell}]]$.

Pour $1\leq i\leq \ell$, on note $\partial_i$ la d\'erivation de
$\C[q_1,\ldots,q_{\ell}]$ d\'efinie par $\partial_i(q_j)=\delta_{ij}$,
$1\leq j\leq \ell$. Elle se prolonge de fa\c con unique en une d\'erivation
de $ \C[[q_1,\ldots,q_{\ell}]]$ qu'on note encore $\partial_i$. On a
alors
$$\Der(\C[q_1,\ldots,q_{\ell}])=\sum_{i=1}^{\ell}\C[q_1,\ldots,q_{\ell}]\partial_i,\quad\mbox{\rm
  et
  }\quad\Der(\C[[q_1,\ldots,q_{\ell}]])=\sum_{i=1}^{\ell}\C[[q_1,\ldots,q_{\ell}]]\partial_i.$$

On note $\g X(\g a/_{W^a})$ (resp. $\g X((\g a/_{W^a}))$) l'ensemble des
d\'erivations de $\C[q_1,\ldots,q_{\ell}]$ (resp. $
\C[[q_1,\ldots,q_{\ell}]]$) qui laissent stable l'id\'eal engendr\'e par
$\it{\Phi}^a|{\g a}$. D'apr\`es la proposition \ref{caspolynomial},  tout \'el\'ement $\delta$ de  $\g X(\g a/_{W^a})$
s'\'ecrit de fa\c con unique $\delta=\sum_{1\leq i\leq \ell}p_i\nabla q_i$,
$p_i\in \C[\g a]^{W^a}$.

On consid\`ere les applications
$$\Der(\C[q_1,\ldots,q_{\ell}]\stackrel{u}{\rightarrow}\C[q_1,\ldots,q_{\ell}]^{\ell}\stackrel{v}{\rightarrow}\C[q_1,\ldots,q_{\ell}]/(\it{\Phi}^a|{\g
  a})\C[q_1,\ldots,q_{\ell}]$$
d\'efinies par $u(\sum_{1\leq i\leq
  \ell}r_i\partial_i)=(r_1,\ldots,r_{\ell})$ et
$v(s_1,\ldots,s_{\ell})=\sum_{1\leq i\leq \ell}s_i\partial_i\it{\Phi}^a|{\g a}$,
de sorte que $\ker v=u(\g X(\g a/_{W^a}))$ ; d'o\`u la suite exacte
\begin{equation}\label{equ8}
0\rightarrow \g X(\g a/_{W^a})\rightarrow
\C[q_1,\ldots,q_{\ell}]^{\ell}\rightarrow
\C[q_1,\ldots,q_{\ell}]/(\it{\Phi}^a|\g a)\C[q_1,\ldots,q_{\ell}].
\end{equation}
 On obtient de m\^eme la suite exacte
%\begin{equation}}\label{equ9}
$$0\rightarrow \g X((\g a/_{W^a}))\rightarrow
\C[[q_1,\ldots,q_{\ell}]]^{\ell}\rightarrow
\C[[q_1,\ldots,q_{\ell}]]/(\it{\Phi}^a|\g a)\C[[q_1,\ldots,q_{\ell}]].$$
%\end{equation}

Comme $\C[[q_1,\ldots,q_{\ell}]]$ est (fid\`element) plat sur
$\C[q_1,\ldots,q_{\ell}]$, on en tire, en tensorisant (\ref{equ8}) par
$\C[[q_1,\ldots,q_{\ell}]]$, que l'application naturelle de $\g X(\g
a/_{W^a})\otimes_{\C[q_1,\ldots,q_{\ell}]}\C[[q_1,\ldots,q_{\ell}]]$
dans $\g X((\g a/_{W^a}))$ est un isomorphisme. La proposition en d\'ecoule.
\end{dem}

\section{Preuve de la premi\`ere partie du th\'eor\`eme
    \ref{principal}}\label{dem1} 

Soit $\mathcal U$ un ouvert compl\`etement $H$-invariant de $\g q$. On
se propose de montrer dans cette section que pour tout $X\in\g
X(\mathcal U)^H$, la d\'erivation $D_X$ de $C^{\infty}(\mathcal U)^H$
laisse stable l'id\'eal $\it{\Phi} C^{\infty}(\g q)^H$. Pour cela il
suffit de montrer que $D_X\it{\Phi}\in \it{\Phi} C^{\infty}(\g q)^H$,
ou, ce qui revient au même, que $\it{\Phi}$ divise $D_X\it{\Phi}$.

Si $\g a\in \Car(\g q)$, on note $\it{\Phi}_{\g a}$ la restriction de
$\it{\Phi}$ \`a $\g a$, qu'on regarde aussi, selon nos conventions, comme  un polyn\^ome sur $\g a_{\C}$.

On note $\g q_{\reg}$ l'ensemble des \'el\'ements semi-simples
r\'eguliers de $\g q$, et si $F$ est une partie de $\g q$, on note
$F_{\reg}=F\cap\g q_{\reg}$. Soit $\g a\in\Car(\g q)$. Comme toute racine de $\g a_{\C}$ dans $\g g_{\C}$ est un multiple d'un \'el\'ement de $\Delta_r(\g g_{\C},\g a_{\C})$ (notation de la section \ref{principal}), il s'ensuit que $\g q_{\reg}$ est l'ensemble des $x\in \g q$ tels que ${\it\Phi}(x)\neq 0$.

\begin{lemme}\label{etape}{Soit $X\in\g X(\mathcal U)^H$. Pour tout
    $\g a\in\Car(\g q)$, il existe $g_{\g a}\in C^{\infty}(\mathcal
    U\cap \g a)$ tel que $(X\it{\Phi}){|(\mathcal U\cap\g a)}=g_{\g a} \it{\Phi}_{\g a}$.}
\end{lemme}

\begin{dem} 
On montre le lemme par r\'ecurrence sur $\dim \g g$. C'est \'evident si
$\g g$ est ab\'elienne. On suppose le
r\'esultat d\'emontr\'e pour les alg\`ebres de Lie r\'eductives  $\g
g'$ telles que $\dim \g g'<\dim \g g$.

On fixe $\g a\in \Car (\g q)$. Pour montrer que $\it{\Phi}_{\g a}$ divise
$(X\it{\Phi}){|\mathcal U\cap \g a}$, il suffit, d'apr\`es le th\'eor\`eme \ref{malgrange}, de montrer que pour tout $a\in\g a$, la s\'erie de Taylor  en $a$ de  $\it{\Phi}_{\g a}$ divise celle de $(X{\it\Phi}){|\mathcal U\cap \g a}$.

Soit $a\in\g a$. On suppose que $a$ n'est pas dans le centre de $\g
g$ et on reprend les notations des sections \ref{descente} et \ref{localisationder}.
L'orthogonal $\g m$ de $\g q^a$ dans $\g
 q$ relativement \`a $\kappa$ est \'egal \`a $[a,\g h]$, l'espace tangent en $a$ de la $H$-orbite de
$a$.

On fixe un voisinage ouvert $\mathcal V$ de $a$ dans $\g q^a$
v\'erifiant les propri\'et\'es de la proposition \ref{desprop}, et on note
$\mathcal W=\pi(H\times \mathcal V)$. On suppose
que  $\det (\ad y)_{\g  g/\g g^a}\neq 0$ et $\Psi^a(y)\neq 0$ pour tout
$y\in\mathcal V$ ; en fait par construction $\mathcal V$ v\'erifie ces
deux propri\'et\'es qui sont clairement \'equivalentes. La
restriction de $X$ \`a $\mathcal V$ s'\'ecrit de fa\c con unique  comme
somme d'un champ de vecteurs $X^a$ sur $\mathcal V$ (\`a valeurs dans
$\g q^a_{\C}$) et d'une
application $Y^a\in C^{\infty}(\mathcal V, \g m_{\C})$. Il est clair que
$X^a$ est $H^a$-invariant, et que, pour tout $x\in
\mathcal V$, l'espace tangent en $x$ de la $H$-orbite de $x$ contient
$\g m$ (voir la preuve du lemme \ref{restcv}), donc $Y^a(x)$ est dans
 le complexifi\'e de l'espace tangent \`a la $H$-orbite de $x$ pour tout
$x\in\mathcal V$. Il s'ensuit que pour toute fonction $f\in
C^{\infty}(\mathcal W)^H$, $(Xf){|\mathcal V}=X^a(f{|\mathcal
  V})$. En particulier, $(X\it{\Phi}){|\mathcal
  V}=X^a(\it{\Phi}|\mathcal V)$.

L'hypoth\`ese de r\'ecurrence implique que $\it{\Phi}^a$ divise $(X^a\it{\Phi}^a)|(\mathcal V\cap\g a)$. On
a
$$(X\it{\Phi}){|\mathcal V}= X^a(\it{\Phi}{|\mathcal V})=
X^a((\it{\Phi}^a\Psi^a)|\mathcal V)=(X^a(\it{\Phi}^a|\mathcal
V)\Psi^a+\it{\Phi}^a(X^a(\Psi^a)\mathcal V),$$
donc, puisque $\g a$ est un sous-espace de Cartan de $\g q^x$, l'hypoth\`ese de r\'ecurrence implique
que $\it{\Phi}^a_{\g a}$ divise $(X\it{\Phi}){|(\mathcal V\cap \g a)}$ ; et comme $\Psi^a$ ne
s'annule pas sur $\mathcal V$, il en d\'ecoule que ${\it \Phi}_{\g a}$
divise $(X\it{\Phi})|(\mathcal V\cap \g a)$. Il en d\'ecoule en
  particulier 
 que  $\it{\Phi}_{\g a}$ divise la  s\'erie de Taylor  en $a$ de
 $(X\it{\Phi}){|(\mathcal U\cap \g a}$. 

On suppose maintenant que $a$ est dans le centre de $\g g$. Alors
$\it{\Phi}_{\g a}$ est un polyn\^ome homog\`ene en $x-a$ de
degr\'e $d$ \'egal au cardinal de  $\Delta_r(\g g_{\C},\g a_{\C})$, et
invariant par $W(\g a)$. 

On note $\g a^{\perp}$ l'orthogonal de $\g a$ dans $\g q$ relativement
\`a $\kappa$. On a $\g
q=\g a\oplus \g a^{\perp}$, car la restriction de $\kappa$
\`a $\g a$ est non d\'eg\'en\'er\'ee. La restriction de $X$ \`a 
$\mathcal U\cap\g a$ s'\'ecrit
alors de fa\c con unique comme somme d'un champ de vecteurs $X_{\g a}$
sur
$\mathcal U\cap\g a$ et d'une application $X^{\perp}\in
C^{\infty}(\mathcal U\cap\g a,\g
a_{\C}^{\perp})$. L'espace tangent en  $x\in\g a_{\reg}$ \`a sa $H$-orbite
s'identifie \`a $\g a^{\perp}$, donc pour tout $f\in
C^{\infty}(\mathcal U)^H$, $Xf(x)=X_{\g a}(f{|\mathcal U\cap \g a})(x)$, et par densit\'e de $\g a_{\reg}$
dans $\g a$, on obtient
$$(Xf){|(\mathcal U\cap\g a}=X_{\g a}(f{|\mathcal U\cap\g a}).$$

Rappelons la d\'erivation $T_aX$ de $\C[[\g q_a]]$
d\'efinie dans la section 3. Vue comme \'el\'ement de $\C[[\g
    q_a]]\otimes_{\C}\g q_{\C}$, c'est une famille $(X_k)_{k\in \N}$
de champs de vecteurs sur $\g q$, chaque $X_k$ est homog\`ene en $x-a$
de degr\'e $k$. Il est clair d'apr\`es leur d\'efinition que les $X_k$ sont
$H^a$-invariants. Cela montre que les composantes homog\`enes en $x-a$
de $T_aX_{\g a}$, qu'on note $X_{\g a,k}$,  sont $W(\g a)$-invariants.

D'apr\`es la proposition
\ref{caspolynomial} (ici ${\it \Phi}^a={\it \Phi}$ et $W^a=W(\g a)$),
pour tout $k\in\N$,  il existe un polyn\^ome $R_k$ tel que
$X_{\g a,k}\it{\Phi}_{\g a}=\it{\Phi}_{\g a}R_k$. Il est clair que ce polyn\^ome est
homog\`ene en $x-a$ de degr\'e $k-1$ ; en particulier il est nul pour
$k=0$. On note $S=(S_k)_{k\in \N}$  l'\'el\'ement de $\C[[\g a_a]]$  tel
que $S_k=R_{k+1}$. Alors $T_a(X_{\g a}\it{\Phi}_{\g a})=\it{\Phi}_{\g a} S$.

On a donc montr\'e que la s\'erie de Taylor de $ 
X_{\g a}\it{\Phi}_{\g a}$ en chaque point $a\in \mathcal U\cap\g a$ est disvisible par
$\it{\Phi}_{\g a}$. Il existe donc, d'apr\`es le th\'eor\`eme \ref{malgrange},
une fonction $g_{\g a}\in C^{\infty}(\mathcal U\cap\g a)$ telle que $X_{\g
    a}\it{\Phi}_{\g a}=\it{\Phi}_{\g a}g_{\g a}$.
\end{dem}

La premi\`ere partie du th\'eor\`eme \ref{principal} est une cons\'equence
imm\'ediate du lemme \ref{etape} et de la proposition \ref{divisioninvariante}.

\section{Preuve de la seconde partie du th\'eor\`eme
    \ref{principal}}\label{dem2}

On fixe un ouvert $\mathcal U$ de $\g q$ compl\`etement
$H$-invariant.

\begin{lemme}\label{restder}{Soit $D\in\Der(C^{\infty}(\mathcal U)^H)$. Pour tout $\g a\in \Car(\g q)$, il  existe une unique d\'erivation $D_{\g a}\in\Der(C^{\infty}(\mathcal  U\cap \g a)^{\inv})$ telle que $D_{\g a}(f{|(\mathcal U\cap \g a)})=(Df){|(\mathcal U\cap \g a)}$
    pour tout $f\in C^{\infty}(\mathcal U)^H$.}
\end{lemme}

\begin{dem} On fixe $\g a\in \Car(\g q)$. L'application $f\mapsto
  f{|(\mathcal U\cap\g a})$ \'etant un morphisme de $\C$-alg\`ebres surjectif  de  $C^{\infty}(\mathcal U)^H$ sur $C^{\infty}(\mathcal U\cap \g a)^{\inv}$, il suffit de
  montrer que $D$ laisse stable son noyau.

Soit $x\in \mathcal U\cap\g a_{\reg}$. D'apr\`es le corollaire \ref{descor}, il
existe 
 $\chi\in C^{\infty}(\mathcal U)^H$ \'egale \`a $1$ dans un  voisinage
de $x$ et de support inclus dans $\pi(H\times \g a_{\reg})$ (notations
de la proposition \ref{desprop}).

Soit $f\in C^{\infty}(\mathcal U)^H$ telle que $f{|(\mathcal U\cap \g a)}$ est nulle. Alors $\chi f$ est nulle. D'o\`u
$$0=D(\chi f)(x)=(D\chi)(x)f(x)+\chi(x)(Df)(x)=(Df)(x).$$
Donc $Df$ est nulle sur $\mathcal U\cap \g a_{\reg}$ ; par suite $Df{|(\mathcal U\cap
  \g a)}$ est nulle, car $\mathcal U\cap \g a_{\reg}$ est dense dans
$\mathcal U\cap\g a$.
\end{dem}

On fixe pour toute la suite de cette section une d\'erivation $D$ de
$C^{\infty}(\mathcal U)^H$ qui laisse stable l'id\'eal $\it{\Phi}
C^{\infty}(\mathcal U)^H$. Donc, d'apr\`es le lemme \ref{restder}, pour tout $\g a\in\Car(\g
q)$, $D_{\g a}\it{\Phi}_{\g a}\in \it{\Phi}_{\g a}C^{\infty}(\mathcal U\cap\g a)^{\inv}$.
 
On fixe $\g a\in\Car(\g q)$ et on note $\mathcal V=\mathcal U\cap \g a
$. On rappelle les g\'en\'erateurs homog\`enes $P_1,\ldots,P_{\ell}$ de
$\C[\g q]^H$. Pour $1\leq i\leq\ell$, on note $p_i$ la
restriction de $P_i$ \`a $\g a$, et $\nabla p_i$ le gradient de $p_i$
relativement \`a $\kappa$. Il d\'ecoule facilement du lemme
\ref{restcv} que $\nabla P_i|\g a=\nabla p_i$.

Comme, pour tout $x\in\g q_{\reg}$,  les vecteurs $\nabla P_1(x),\ldots, \nabla P_{\ell}(x)$ sont lin\'eairement
ind\'ependants (voir \cite{K-R}, Th\'eor\`eme 13), il en est donc de m\^eme des
vecteurs $\nabla p_1(x),\ldots, \nabla p_{\ell}(x)$, $x\in \g a_{\reg}$. Soit
$x\in \mathcal V_{\reg}$. On a $\g q^x=\g a$ et $H^x$ op\`ere trivialement dans
$\g a$. Il existe donc d'apr\`es la proposition \ref{desprop} un
voisinage ouvert $\mathcal O$ de $x$ dans $\g a$, inclus dans
$\mathcal V_{\reg}$, tel que $\mathcal W=\pi(H\times \mathcal O)$ soit un
ouvert totalement $H$-invariant de $\g q$, inclus dans $\mathcal U$,
et l'application de restriction de $C^{\infty}(\mathcal W)^H$ dans
$C^{\infty}(\mathcal O)$ soit un isomorphisme. La restriction de $D$ \`a
$C^{\infty}(\mathcal W)^H$ (proposition \ref{locder}) d\'efinit une d\'erivation
de $C^{\infty}(\mathcal O)$ ; cette d\'erivation s'\'ecrit de fa\c con
unique sous la forme 
$\sum_{1\leq i\leq\ell}f_i\nabla p_i$, $f_i\in C^{\infty}(\mathcal
O)$, car les vecteurs $\nabla p_1(y),\ldots,\nabla p_{\ell}(y)$
forment une base de $\g a$ pour tout $y\in\mathcal O$. On en d\'eduit
qu'il existe  des fonctions
$\varphi_1,\ldots,\varphi_{\ell}\in C^{\infty}(\mathcal V_{\reg})$
uniquement d\'etermin\'ees telles que 
$$(D_{\g a}g){|\mathcal V_{\reg}}= \sum_{i=1}^{\ell}\varphi_i\nabla
  p_i(g{|\mathcal V_{\reg}}),\quad\mbox{\rm pour tout } g\in
  C^{\infty}(\mathcal V)^{\inv}.$$

Appliquant ceci aux fonctions $p_j$, $1\leq j\leq \ell$, en chaque
$x\in\mathcal V_{\reg}$, on trouve le syst\`eme lin\'eaire
\begin{equation*}\label{systeme}
D_{\g a}p_j(x)=\sum_{i=1}^{\ell}\varphi_i(x)\nabla p_i\cdot p_j(x).
\end{equation*}

Il existe une constante $c\neq 0$ telle que 
$$\det (\nabla p_i\cdot p_j(x))_{1\leq i,j\leq\ell}=c\it{\Phi}_{\g a}(x),\quad\mbox{\rm pour tout } x\in\g a\; ; $$
cela d\'ecoule des propri\'et\'es des polyn\^omes invariants par un groupe
 engendr\'e par des r\'eflexions (\cite{Bou}, chapitre V, §5,
  proposition 5). 

Les formules de Cramer montrent qu'il existe des fonctions
  $\psi_1,\ldots,\psi_{\ell}\in C^{\infty}(\mathcal V)^{\inv}$ telles que,
  pour tout $1\leq i\leq \ell$, 
$$\it{\Phi}_{\g a}(x)\varphi_i(x)=\psi_i(x),\quad\mbox{\rm pour tout }
  x\in\mathcal V_{\reg}.$$
Il s'ensuit que
\begin{equation}\label{PhiD}
\it{\Phi}_{\g a}D_{\g a}(f)=\sum_{1\leq i\leq\ell}\psi_i\nabla p_i(f),\quad
\mbox{\rm pour tout  }f\in C^{\infty}(\mathcal V)^{\inv}.
\end{equation}

Remarquons que $\nabla p_i\cdot p_j(x)=\nabla P_i\cdot P_j(x)$ pour tout $x\in\g
a$, il s'ensuit que les polyn\^omes $H$-invariants $\det (\nabla P_i\cdot P_j)_{1\leq
  i,j\leq\ell}$ et $c\it{\Phi}$ co\"incident sur $\g a$, ils  sont donc \'egaux. Les formules
de Cramer montrent alors que  les $\psi_i$ sont des restrictions \`a
$\mathcal V$
de fonctions $\widetilde{\psi}_i\in C^{\infty}(\mathcal U)^H$. Plus
pr\'ecis\'ement, les $\widetilde{\psi}_i$ sont uniquement d\'etermin\'ees par
$$\it{\Phi}DP_j=\sum_{1\leq i\leq\ell}\widetilde{\psi}_i\nabla
P_i\cdot P_j,\quad 1\leq j\leq \ell.$$

On veut montrer que les fonctions $\varphi_i$ se prolongent en des
fonctions de classe $C^{\infty}$ sur $\mathcal V$ ou, ce qui revient
au m\^eme, que les $\psi_i$ sont divisibles par $\it{\Phi}_{\g a}$. Pour cela, en vertu du
th\'eor\`eme \ref{malgrange},  il suffit de
montrer que les  s\'eries de Taylor des $\psi_i$ en chaque point $a\in
\mathcal V$ sont divisibles par $\it{\Phi}_{\g a}$.

On fixe $a\in \mathcal V$ et on reprend les notations de la section
 \ref{localisationder}. En particulier $(q_1,\ldots,q_{\ell})$ est un
 ensemble de g\'en\'erateurs homog\`enes en $x-a$ de $\C[\g a]^{W^a}$.

\begin{lemme}\label{fin}{Il existe des  polyn\^omes $m_{ij},
 1\leq i,j\leq \ell$, dans $\C[\g a]^{W^a}$ tels que, pour tout $1\leq j\leq \ell$, on ait
$$\nabla p_j=\sum_{1\leq i\leq \ell}m_{ij}\nabla q_i.$$
De plus la matrice $(m_{ij})_{1\leq i,j\leq \ell}$ est inversible dans
$M(\ell,\C[[\g a_a]]^{W^a})$.}
\end{lemme}
\begin{dem} Comme les $\nabla p_i$ sont $W$-invariants et donc
  $W^a$-invariants, l'existence des $m_{ij}$ d\'ecoule de \cite{So}.

On a d\'ej\`a vu au cours de la preuve du lemme \ref{taylorinv} qu'il existe des
  polyn\^omes $Q_1,\ldots,Q_{\ell}\in\C[\g q^a]^{H^a}$ tels que
  $Q_i|\g a=q_i$. Il d\'ecoule facilement du lemme \ref{restcv} que, pour tout $1\leq i\leq
  \ell$, $\nabla Q_i|\g a=\nabla q_i$. On note $M_{ij}$ l'unique
  polyn\^ome $H^a$-invariant sur $\g q^a$ tel que $M_{ij}|\g
  a=m_{ij}$. Les restrictions des $\nabla P_i$ \`a $\g q^a$ sont \`a
  valeurs dans $\g q^a_{\C}$ d'apr\`es le lemme \ref{restcv} ; comme ces
  restrictions sont $H^a$-invariantes, il d\'ecoule de la d\'efinition des
  $M_{ij}$ que, pour tout $1\leq j\leq \ell$,
 $$\nabla P_j|\g q^a=\sum_{1\leq i\leq \ell}M_{ij}\nabla Q_i.$$

On note $\g r$ l'ensembles des \'el\'ements r\'eguliers de $\g q_{\C}$ (ce
sont les \'el\'ements dont l'orbite sous $H_{\C}$ est de dimension maximale ; ils
ne sont pas forc\'ement semi-simples). D'apr\`es (\cite{K-R},
th\'eor\`eme 13), pour tout $x\in \g r$, les  vecteurs $\nabla
P_1(x),\ldots,\nabla P_{\ell}(x)$ sont lin\'eairement ind\'ependants. De
m\^eme si l'on note $\g r^a$ l'ensemble des \'el\'ements r\'eguliers de $\g
q_{\C}^ a$, c'est-\`a-dire l'ensemble des \'el\'ements dont l'orbite sous
$H_{\C}^a$ est de dimension maximale, alors pour tout $y\in\g r^a$, les vecteurs $\nabla Q_1(y),\ldots,\nabla Q_{\ell}(y)$ sont
lin\'eairement ind\'ependants ; comme $\g r \cap \g q_{\C}^a\subset \g
r^a$, cette propri\'et\'e est en particulier vraie pour tout $y\in \g
r\cap \g q_{\C}^a$. Il existe $n\in\g q_{\C}^a$ nilpotent tel que
$x=a+n\in \g r\cap \g q_{\C}^a$ (voir le d\'ebut de la preuve du
th\'eor\`eme 8 de \cite{K-R}). Alors la matrice $(M_{ij}(x))$ est
inversible. Mais comme ses \'el\'ements sont $H^a$-invariants, on a
$M_{ij}(x)=M_{ij}(a)$ ; donc $\det (M_{ij})(a)\neq 0$, cela prouve
que $\det (m_{ij})$ est inversible dans $\C[[\g a]]^{W^a}$, et donc la
matrice 
$(m_{ij})$ est inversible dans $ M(\ell,\C[[\g a]]^{W^a})$.
\end{dem}

D'apr\`es le lemme \ref{idealpreserv} et la
 proposition \ref{derformel},  la d\'erivation $(D_{\g a})_a$ de $\C[[\g
 a_a]]^{W^a}$ est de la forme $\sum_{1\leq i\leq \ell}S_i\nabla q_i$,
 $S_i\in \C[[\g a_a]]^{W^a}$. 

On consid\`ere la d\'erivation $\delta=\it{\Phi}_{\g a}D_{\g a}$. Alors 
$$\delta_a=\sum_{1\leq i\leq \ell}\it{\Phi}_{\g a}S_i\nabla q_i,$$
et, d'apr\`es (\ref{PhiD}) et le lemme \ref{fin},
\begin{eqnarray*}
\delta_a&=&\sum_{1\leq j\leq \ell}(T_a\psi_j)\nabla p_j\\
&=&\sum_{1\leq j\leq \ell}(\sum_{1\leq i\leq \ell}(T_a\psi_j)m_{ij}\nabla q_i)\\
&=&\sum_{1\leq i\leq \ell}(\sum_{1\leq j\leq
  \ell}m_{ij}(T_a\psi_j))\nabla q_i.
\end{eqnarray*}
Donc par unicit\'e de l'\'ecriture de $\delta_a$ (proposition \ref{derformel}), pour
tout $1\leq i\leq \ell$, on a
$$\it{\Phi}_{\g a}S_i=\sum_{1\leq j\leq
  \ell}m_{ij}(T_a\psi_j).$$
Comme la matrice $(m_{ij})$ est inversible dans $M(\ell,\C[[\g a_a]]^{W^a})$,
  il en est de m\^eme de sa transpos\'ee. Les formules de Cramer
  montrent alors que $T_a\psi_j$ est divisible par $\it{\Phi}_{\g a}$ dans
  $\C[[\g a_a]]^{W^a}$.

Ce qui pr\'ec\`ede montre  que, pour tout $1\leq i\leq \ell$ et
pour tout $\g a\in\Car(\g q)$, la
restriction de $\widetilde{\psi_i}$ \`a $\mathcal U\cap\g a$ est divisible par $\it{\Phi}_{\g
  a}$. Il existe donc, d'apr\`es la proposition
\ref{divisioninvariante}, des fonctions 
$\widetilde{\varphi}_1,\ldots,\widetilde{\varphi}_{\ell}\in
C^{\infty}(\mathcal U)^H$ telles que
$\widetilde{\psi}_i=\it{\Phi}\widetilde{\varphi}_i$. On pose $X=\sum_{1\leq i\leq
  \ell}\widetilde{\varphi}_i\nabla P_i$.

Soit $f\in C^{\infty}(\mathcal U)^H$. Pour tout $\g a\in \Car(\g q)$,
les restrictions \`a $\mathcal U\cap \g a_{\reg}$ des fonctions  $Df$ et  $Xf$
co\"incident ; il s'ensuit que $Df$ et $Xf$ co\"incident sur $\mathcal
U_{\reg}$, et donc $Df=Xf$. D'o\`u $D_X=D$. L'unicit\'e des
$\widetilde{\varphi}_i$ est claire d'apr\`es la construction.

\section{Remarques}

\begin{remarque} { Soient $P_1,\ldots,\P_{\ell}$, $\ell$ le rang de
    $\g g$, un syst\`eme minimal de g\'en\'erateurs homog\`enes de $\C[\g
    g]^G$. Il d\'ecoule par exemple de (\cite{B}, corollaire 3.5) que si
    $\mathcal U$ est un ouvert compl\`etement $G$-invariant de $\g g$, alors tout
    \'el\'ement de $\g X(\mathcal U)^G$ s'\'ecrit de fa\c con unique 
   sous la forme $\sum_{1\leq i\leq\ell}f_i\nabla P_i$, avec $f_i\in
    C^{\infty}(\mathcal U)^G$. On note $\it{\Phi}=\det (\nabla P_i\cdot
    P_j)_{1\leq i,j\leq\ell}$. Le th\'eor\`eme \ref{principal}, appliqu\'e
    au cas diagonal, se traduit ainsi :  l'application $X\mapsto D_X$ de $\g
    X(\mathcal U)^G$ dans $\Der(C^{\infty}(U)^G$ est un isomorphisme
    de $\g  X(\mathcal U)^G$ sur le sous-espace de
    $\Der(C^{\infty}(U)^G$ form\'e des d\'erivations qui pr\'eservent
    l'id\'eal engendr\'e par $\it{\Phi}$.}
\end{remarque}

\begin{remarque} {Soit $\mathcal U$ un ouvert compl\`etement
    $H$-invariant de $\g q$. On note $\mathscr Z(\mathcal U)$
    l'ensemble des \'el\'ements de $\g X(\mathcal U)^H$ qui annulent
    tous les \'el\'ements de $C^{\infty}(\mathcal U)^H$. C'est un
    sous-$C^{\infty}(\mathcal U)^H$-module de $\g X(\mathcal U)^H$, et
    il d\'ecoule du
    th\'eor\`eme \ref{principal} que $\mathscr Z(\mathcal U)$ est un
    suppl\'ementaire du sous-module $\sum_{1\leq i\leq
    \ell}C^{\infty}(\mathcal U)^H\nabla P_i$. On peut voir aussi que
    si $X\in\g X(\mathcal U)^H$, alors $X\in\mathscr Z(\mathcal U)$ si
    et seulement si pour tout $\g a\in\Car(\g q)$
    et pout tout $x\in \g a$, $X(x)$ appartient \`a l'orthogonal de $\g
    a_{\C}$ dans $\g q_{\C}$ relativement \`a $\kappa$. L'exemple suivant montre
    qu'en g\'en\'eral $\mathscr Z
    (\mathcal U)$ n'est pas engendr\'e par des champs de vecteurs
    invariants \`a coefficients polynomiaux. 
}
\end{remarque}

\begin{exemple}{Soient $G=\SL(3,\R)$, $\g g=\g{sl}(3,\R)$ son alg\`ebre de
 Lie et $\sigma$ l'involution de $G$ d\'efinie par
$$\sigma(M)=I_{2,1}\ ^{t}M^{-1}I_{2,1}, \quad M\in \SL(3,\R),\quad
 \mbox{\rm{o\`u}  }\quad I_{2,1}=\diag(1,1,-1).$$
Le groupe des points fixes de $\sigma$ est $SO(2,1)$ ; on note $H$ sa
composante neutre.

La diff\'erentielle de $\sigma$, qu'on note encore $\sigma$, est d\'efinie
 par
$$\sigma(A)=-I_{2,1}\ ^{t}AI_{2,1}, \quad A\in\g{sl}(3,\R).$$
Avec les notations de la section \ref{resultat}, on a
$$\g
h=\left\{\left(\begin{array}{ccc}0&b&c\\-b&0&d\\c&d&0\end{array}\right)
    ;b,c,d\in \R\right\}, \g q=\left\{\left(\begin{array}{ccc}a&b&c\\b&d&e\\-c&-e&-(a+d)\end{array}\right)
    ;a,b,c,d,e\in \R\right\}.$$

Soit 
$$\g
a=\left\{\left(\begin{array}{ccc}x&0&y\\0&-2x&0\\-y&0&x\end{array}\right);
  x,y\in\R\right\} ;$$
c'est un sous-espace de Cartan de $\g q$.

On note avec un indice $\C$ les complexifi\'es des objets d\'efinis ci-dessus ; ainsi
$G_{\C}=\SL(3,\C)$, $H_{\C}$ l'ensemble des points fixes de
$\sigma_{\C}$ ; $H_{\C}$ est isomorphe \`a $SO(3,\C)$, donc
connexe. On note $M_{\C}$ le centralisateur de $\g a_{\C}$ dans
$H_{\C}$. On a
$$M_{\C}=\left\{\left(\begin{array}{ccc}1&0&0\\0&1&0\\0&0&1\end{array}\right)
  ; \left(\begin{array}{ccc}-1&0&0\\0&1&0\\0&0&-1\end{array}\right)
    ;\left(\begin{array}{ccc}0&0&i\\0&-1&0\\-i&0&0\end{array}\right)
      ;\left(\begin{array}{ccc}0&0&-i\\0&-1&0\\i&0&0\end{array}\right)\right\}.$$D'o\`u $\g q_{\C}^{M_{\C}}=\g a_{\C}$.

D'apr\`es \cite{K-R}, l'ensemble $\g X_p(\g q_{\C})^{H_{\C}}$ des champs de vecteurs \`a coefficients
polynomiaux et $H_{\C}$-invariants sur $\g q_{\C}$ est un $\C[\g
  q_{\C}]^{H_{\C}}$-module libre de rang $\dim \g
q_{\C}^{H_{\C}}=2$. Pour \^etre pr\'ecis, dans \cite{K-R} les auteurs
consid\`erent l'action dans $\g q_{\C}$ du groupe $K_{\C}$, commutant de
$\sigma_{\C}$ dans le groupe adjoint de $\g g_{\C}$ , mais cela n'a
pas d'importance ici car le groupe $K_{\C}$ est connexe, et donc les
champs de vecteurs polynomiaux invariants sous l'action de $H_{\C}$ ou
sous l'action de $K_{\C}$ sont les m\^emes. On en d\'eduit que $\g
X_p(\g q_{\C})^{H_{\C}}=\C[\g q_{\C}]^{H_{\C}}\nabla P_1+\C[\g
  q_{\C}]^{H_{\C}}\nabla P_2$, o\`u $(P_1,P_2)$ est un syst\`eme de
g\'en\'erateurs homog\`enes de $\C[\g q_{\C}]^{H_{\C}}$ ; ceci d\'ecoule aussi
de \cite{SB}, car la paire $(\g g_{\C},\g h_{\C})$
peut \^etre vue comme la complexfi\'ee d'une paire riemannienne
sym\'etrique non hermitienne (la repr\'esentation de $\g h_{\C}$ dans $\g
q_{\C}$ est irr\'eductible). En particulier, on a $\g
X_p(\g q_{\C})^{H_{\C}}\cap \mathscr Z(\g q)=\{0\}$.

Le centralisateur $M$ de $\g a$ dans $H$ est \'egal \`a 
$$\left\{\left(\begin{array}{ccc}1&0&0\\0&1&0\\0&0&1\end{array}\right)  ;
 \left(\begin{array}{ccc}-1&0&0\\0&1&0\\0&0&-1\end{array}\right)\right\},$$
donc 
$$\g
q^{M}=\left\{\left(\begin{array}{ccc}x&0&y\\0&z&0\\-y&0&-(x+z)\end{array}\right)
;x,y,z\in \R\right\}.$$
On note $v$ l'\'el\'ement de $\g q^M$ correspondant \`a $x=1$ et $y=z=0$ ;
il appartient \`a l'orthogonal de $\g a$ dans $\g q$ relativement \`a $\kappa$.

Soit $x\in\g a_{\reg}$ et soit $\mathcal V$ un voisinage ouvert de $x$
dans $\g a_{\reg}$ comme dans la proposition \ref{desprop}. Soit $\varphi\in
C^{\infty}(\mathcal V)$. L'application $y\mapsto \varphi(y)v$ se
prolonge de fa\c con unique en un \'el\'ement $X$ de $\g X(\mathcal
W)^H$. Ce champ de vecteurs annule toutes les fonctions dans
$C^{\infty}(\mathcal W)^H$ et n'appartient pas au sous-$C^{\infty}(\mathcal U)^H$-module de $\g
X(\mathcal U)^H$ engendr\'e par les champs de vecteurs invariants \`a
coefficients polynomiaux.}
\end{exemple}

\end{document}